\RequirePackage{lineno}
\documentclass[aps,titlepage,11pt,superscriptaddress]{revtex4}
\usepackage{epsfig}
\usepackage{float}
\usepackage{graphicx}
\usepackage{times}
\usepackage{amsmath}
\usepackage{color}
\usepackage{lineno}

\begin{document}


\title{The paradigm of tax-reward and tax-punishment strategies in the advancement of public resource management dynamics}

\author{Lichen Wang}
\affiliation{College of Science, Northwest A \& F University, Yangling 712100, China}

\author{Yuyuan Liu}
\affiliation{College of Science, Northwest A \& F University, Yangling 712100, China}

\author{Ruqiang Guo}
\affiliation{College of Science, Northwest A \& F University, Yangling 712100, China}

\author{Liang Zhang}
\affiliation{College of Science, Northwest A \& F University, Yangling 712100, China}

\author{Linjie Liu}
\email{linjie140@126.com}
\affiliation{College of Science, Northwest A \& F University, Yangling 712100, China}
\affiliation{College of Economics \& Management, Northwest A\&F University, Yangling 712100, China}

\author{Shijia Hua}
\email{sjhua@nwafu.edu.cn}
\affiliation{College of Science, Northwest A \& F University, Yangling 712100, China}

\begin{abstract}\noindent
\\
In contemporary society, the effective utilization of public resources remains a subject of significant concern. A common issue arises from defectors seeking to obtain an excessive share of these resources for personal gain, potentially leading to resource depletion. To mitigate this tragedy and ensure sustainable development of resources, implementing mechanisms to either reward those who adhere to distribution rules or penalize those who do not, appears advantageous. We introduce two models: a tax-reward model and a tax-punishment model, to address this issue. Our analysis reveals that in the tax-reward model, the evolutionary trajectory of the system is influenced not only by the tax revenue collected but also by the natural growth rate of the resources. Conversely, the tax-punishment model exhibits distinct characteristics when compared to the tax-reward model, notably the potential for bistability. In such scenarios, the selection of initial conditions is critical, as it can determine the system's path. Furthermore, our study identifies instances where the system lacks stable points, exemplified by a limit cycle phenomenon, underscoring the complexity and dynamism inherent in managing public resources using these models.
\end{abstract}

\maketitle

\noindent

The rampant overexploitation of public resource pools has emerged as a pervasive issue in today's world, cutting across geographical and cultural boundaries \cite{plus1,plus2,plus3,gg1}. As humanity's footprint on ecosystems escalates, the reciprocal influence of ecological variables on human societies is concurrently intensifying. This intricate interplay between societal and ecological systems is not an isolated incident but a globally prevalent phenomenon \cite{liu2023coevolutionary,rl2,third2,rll1}. A case in point is the widespread overgrazing of public rangelands. Across numerous regions, rangeland resources are subjected to relentless exploitation, culminating in their alarming depletion \cite{gdfm1,gdfm2}. This predicament frequently arises in the absence of effective management and vigilant supervision, as herders, driven by their necessities, overgraze the land, thereby impeding the resources' ability to rejuvenate, and propelling them towards eventual exhaustion. This scenario is not confined to rangelands alone but is mirrored in a multitude of contexts, such as the overfishing of marine resources, the excessive consumption of water resources, the exacerbation of global climate change \cite{ref21,ref19,ref16,ref17,wang2020communicating}. These instances underscore the urgent need for comprehensive strategies to manage public resources sustainably.

To address the challenge of resource overexploitation, managers often establish allocation rules. Despite these efforts, individual desires for greater benefits frequently lead to the violation of these rules, exacerbating resource depletion. This scenario underscores the significance of cooperative behavior, wherein participants comply with allocation rules out of a collective commitment to resource sustainability. This raises a pivotal question: How can we motivate individuals to persist in such cooperative behavior for the enduring sustainability of resources? This inquiry has generated substantial interest and debate in the field \cite{hz1,hz2,hz3,capraro2021mathematical}. 

The intricacies of public resource management often call for the adoption of reward and punishment strategies, a stance substantiated by a wealth of studies \cite{sigmund2010social,hilbe2010incentives,ref4,ref5,ref7,ref8,ref10,ref12,ref13,jl1,wang2023emergence,jl2}. In particular, the work of Chen and Szolnoki highlights the effectiveness of punitive mechanisms \cite{chenpcb2018}. Their inspection-punishment model, which targets and penalizes defectors, demonstrates how different natural growth rates impact the evolution of public resource pools. Yet, this model applies a constant severity of punishment, while in reality, levels of punishment often adjust in response to varying violation rates. The dynamic incentive allocation mechanism can address these shortcomings. Under this scheme, all individuals pay taxes. In the reward scenario, the total tax collected is evenly distributed among all cooperators as a reward \cite{ref23,hua2024coevolutionary}. In the punishment scenario, the total tax is evenly allocated to the defectors as a penalty \cite{ref24,sun2021combination}. However, the use of this dynamic incentive allocation mechanism is seldom mentioned or employed in research on public pool management. Therefore, its effectiveness in promoting the emergence of cooperation and ensuring sustainable resource use remains to be seen.

This study pioneers the application of a tax incentive strategy to the management of public resource pools, utilizing replicator equations to portray the evolutionary dynamics of cooperation \cite{ref25,ref26,ref27,ref28,ref35}. Concurrently, a modified logistic model is employed to depict the evolutionary dynamics of public resource pools \cite{ref35,ref36}. We construct tax-reward and tax-punishment models to offer a more comprehensive understanding and prediction of the temporal evolution of public resource pools. We take into account the finite and renewable nature of the resource pool, along with the reciprocal influence of individual behavior on resources. For the tax-reward strategy, we demonstrate the temporal evolution of the proportion of cooperators and track the changes in the common resource pool, capturing the dynamics associated with varying tax revenues and natural growth rates. Our findings reveal that the system may ultimately gravitate towards one of three states: full cooperation, full defection, or a coexistence of cooperation and defection. For the tax-punishment strategy, we similarly demonstrate the temporal evolution of the proportion of cooperators and track the changes in the common resource pool, capturing the dynamics associated with varying tax revenues and natural growth rates. Contrary to the tax-reward strategy, we observe that the system may exhibit bistability, suggesting that the stability of the system is contingent upon initial conditions. Furthermore, we identify instances where stable fixed points do not exist, such as the limit cycle discussed in this paper, representing a unique case of unstable fixed points. Lastly, we corroborate the theoretical accuracy and the efficacy of our model through numerical examples and Monte Carlo simulations.

\vbox{}

\leftline{\textbf{Model and Methods}}

Let's envision a group of  $N$ individuals, all utilizing a shared, finite, and renewable resource pool, with a natural growth rate denoted as $r$. The current resource level of the common pool is represented by $y$. Typically, the logistic model is used to describe the dynamic changes \cite{chenpcb2018}. The equation for the current resource total $y$ with respect to time $t$ is as follows: \begin{eqnarray}
	\dot{y} = ry \left(1 - \frac{y}{R_{m}}\right),
\end{eqnarray}
 where $\dot{y} = \dfrac{dy}{dt}$, and $R_{m}$ represents the maximum capacity of the resource environment. To ensure the sustainable development of the common resource, a distribution rule is proposed in reference \cite{chenpcb2018}. Each individual obtains a resource quantity of $b_{l} = \frac{b_{m} y}{R_{m}}$ from the current resource pool, where $b_{m}$ specifies the maximum resource amount an individual is allowed to use when the resource reaches its maximum capacity.

To streamline our model, we posit that individuals adopt one of two possible strategies. The first strategy, undertaken by law-abiders or cooperators, involves adhering to established allocation rules. Conversely, the second strategy is employed by law-breakers or defectors, who disregard these rules in favor of personal gain, securing a portion, $b_v$, that surpasses the designated quota, $b_l$. We define $b_v$ as $b_l (1 + \alpha)$, with $\alpha$ signifying the ``defection rate", akin to the dilemma strength often mentioned in prior research on $2\times2$ and $N$-player games \cite{arefin2020social,wang2015universal,ito2018scaling,tanimoto2007relationship}. In this scenario, we presume that an individual's payoff is equivalent to the share of resources they extract from the collective resource pool.

In analyzing a well-mixed population, we employ replicator equations to model the temporal evolution of competing strategies. The fraction of cooperators within the population is denoted by $x$; consequently, we can assert the governing equation as follows:
\begin{eqnarray}
	\dot{x} = x(1-x)\left(P_{C}-P_{D}\right),
\end{eqnarray} 
wherein $\dot{x}$ signifies the rate of change of $x$ concerning time $t$, symbols $P_C$ and $P_D$ pertain to the average payoffs for cooperators and defectors, correspondingly.

As the concept of coevolution accounts for the reciprocal influence of individual behavior on the prevailing resource conditions, the governing equation governing the abundance of shared resources can be extended to:

\begin{eqnarray}
	\dot{y} = r y\left(1-\frac{y}{R_{m}}\right) - N\left[b_{l} x + (1-x) b_{v}\right].
\end{eqnarray}

To avoid cooperators turning into defectors due to unfair distribution of benefits, leading to the depletion of common resources, we propose the following two incentive strategies to promote the sustainable development of common resources.

\leftline{\textbf{Tax-reward}}
\quad \ In the context of real-world resource management, rewarding cooperators is generally considered an effective strategy. Therefore, we propose a tax-reward strategy. In this strategy, we assume a population consisting of $N$ individuals, with $N_C$ being the number of cooperators. Each individual is required to pay a tax of $\delta$. The total revenue generated from these taxes is then used to provide rewards to the cooperators. The average payoff for cooperators, denoted as $P_C$, and the average payoff for defectors, denoted as $P_D$, can be expressed as follows:
\begin{align}
	\begin{split}
		P_C &= \sum_{k=0}^{N-1}\binom{N-1}{k}x^k(1-x)^{N-1-k}\left(b_l - \delta + \frac{N\delta}{k+1}\right), \\
		P_D &= \sum_{k=0}^{N-1}\binom{N-1}{k}x^k(1-x)^{N-1-k}\left[b_l(1+\alpha) - \delta\right].
	\end{split}
\end{align}
According to the property of binomial coefficients $\binom{N-1}{k}\frac{N}{k+1}=\binom{N}{k+1}$, we can obtain:
\begin{align}
	\begin{split}
		P_C - P_D &= \sum_{k=0}^{N-1}\binom{N-1}{k}x^k(1-x)^{N-1-k}\left(-b_l\alpha + \frac{N\delta}{k+1}\right) \\
		&= -\frac{b_my}{R_m}\alpha + \delta\sum_{k=0}^{N-1}\binom{N}{k+1}x^k(1-x)^{N-1-k} \\
		&= -\frac{b_my}{R_m}\alpha + \delta\sum_{k'=1}^{N-1}\binom{N}{k'}x^{k'-1}(1-x)^{N-k'} \\
		&= -\frac{b_my}{R_m}\alpha + \delta\sum_{k'=1}^{N-1}\binom{N}{k'}\frac{x^{k'}(1-x)^{N-k'}}{x} \\
		&= -\frac{b_my}{R_m}\alpha + \delta\frac{1-(1-x)^N}{x}.
	\end{split}
\end{align}

Based on recent dynamic feedback modeling approaches \cite{chenpcb2018,arefin2021imitation}, we derived our tax-reward model as follows: 
\begin{eqnarray}\label{eq_1}
	\left\{\begin{array}{l}
		\dot{x}=x(1-x)\left[\delta\frac{1-\left(1-x\right)^N}{x}-\frac{y}{R_{m}} b_{m} \alpha\right], \\
		\dot{y}=r y\left(1-\frac{y}{R_{m}} \right)-N \frac{y}{R_{m}} b_{m}[1+(1-x) \alpha] .
	\end{array}\right.
\end{eqnarray}
\leftline{\textbf{Tax-punishment}}
\quad \ In real-life situations, when the proportion of cooperators is relatively high, the cost of rewarding them can often be substantial. In such cases, it becomes necessary to consider the punishment of defectors. We propose a tax-punishment strategy. Firstly, let $N_D$ denote the number of defectors among the $N$ individuals in the population. Each individual is required to pay a tax of $\delta$. The total revenue generated from these taxes is then utilized to administer punishment to the defectors. Similarly, the average payoff for cooperators, denoted as $P_C$, and the average payoff for defectors, denoted as $P_D$, can be expressed as follows:
\begin{align*}
	\begin{split}
	P_C &= \sum_{k=0}^{N-1}\binom{N-1}{k}x^k(1-x)^{N-1-k}\left(b_l - \delta\right), \\
	P_D &= \sum_{k=0}^{N-1}\binom{N-1}{k}x^k(1-x)^{N-1-k}\left[b_l(1+\alpha) - \delta-\frac{N\delta}{N_D+1}\right].
	\end{split}
\end{align*}

According to the property of binomial coefficients $\binom{N-1}{k}\frac{N}{N-k}=\binom{N}{k}$, we can obtain:
\begin{align*}
	\begin{split}
			P_C - P_D &= \sum_{k=0}^{N-1}\binom{N-1}{k}x^k(1-x)^{N-1-k}\left(-b_l\alpha + \frac{N\delta}{N-k}\right) \\
		&= -\frac{b_my}{R_m}\alpha + \delta\sum_{k=0}^{N-1}\binom{N}{k}x^k(1-x)^{N-1-k} \\
		&= -\frac{b_my}{R_m}\alpha + \delta\sum_{k=0}^{N-1}\binom{N}{k}\frac{x^{k}(1-x)^{N-k}}{1-x} \\
		&= -\frac{b_my}{R_m}\alpha + \delta\frac{1-x^N}{1-x}.
	\end{split}
\end{align*}

Based on the above, we obtained our tax-punishment model as follows:
\begin{eqnarray}\label{eq_2}
	\left\{\begin{array}{l}
		\dot{x}=x(1-x)\left(\delta\frac{1-x^N}{1-x}-\frac{y}{R_{m}} b_{m} \alpha\right), \\
		\dot{y}=r y\left(1-\frac{y}{R_{m}} \right)-N \frac{y}{R_{m}} b_{m}[1+(1-x) \alpha] .
	\end{array}\right.
\end{eqnarray}

To help readers easily understand all the parameters and variables introduced in our work, we present them in Table \ref{table1}.
\vspace*{-7pt}
\begin{table}[!h]
	\caption{Symbols and meanings used in this article.}
	\label{table1}
	\begin{tabular}{p{2cm}p{11cm}} 
		\hline
		\textbf{Symbol} & \textbf{Meaning} \\
		\hline
		$N$ & Number of individuals in the game group \\
		$N_C$ & Number of cooperators \\
		$N_D$ & Number of defectors \\
		$y$ & Current environmental resource level \\
		$x$ & Current proportion of cooperators \\
		$\alpha$ & Defection rate of defectors \\
		$\delta$ & Tax paid by each participant \\
		$R_m$ & Maximum resource carrying capacity of the environment \\
		$b_m$ & Maximum amount of resources allowed for an individual when resources reach maximum capacity \\
		$b_l$ & Amount of resources obtained by cooperators \\
		$b_v$ & Amount of resources obtained by defectors \\
		$P_C$ & Average payoff of cooperators \\
		$P_D$ & Average payoff of defectors \\
		\hline
	\end{tabular}
	\vspace*{-4pt}
\end{table}

\vbox{}

\leftline{\textbf{Results}}
\leftline{\textbf{Tax-reward}}

For the tax-reward model, utilizing equation \eqref{eq_1}, there are four boundary fixed points in the system: $\left(x,y\right)=\left(0,0\right)$, $\left(1,0\right)$, $\left(1,R_m-\frac{Nb_m}{r}\right)$, and $\left(0,R_m-\frac{Nb_m}{r}\left(1+\alpha\right)\right)$. Additionally, when the conditions $\delta-\alpha b_m + \frac{\alpha N b_m^2}{r R_m}<0$ and $N\delta-\alpha b_m+\frac{\alpha N b_m^2}{r R_m}+\frac{\alpha^2 N b_m^2}{r R_m}>0$ are satisfied, there exists one and only one interior fixed point $\left(x^\ast,y^\ast\right)$, where $0<x^\ast<1$ and $0<y^\ast<R_m$. otherwise, the system does not have an interior fixed point. We use the sign of the eigenvalues of the Jacobian matrix to determine the stability of the fixed points (detailed proof can be found in Supplementary Information). Now we present the evolution of cooperation proportion and resource abundance over time for different initial conditions, along with their corresponding phase diagrams where the solid black dots represent stable fixed points , while the hollow black dots represent unstable fixed points. Here, we define $E_C=\frac{b_mN}{Rm}$ and $E_D=\frac{b_mN(1+\alpha)}{R_m}$. We will discuss the following three cases based on different natural growth rates $r$.

\begin{figure}[!h]
	\centering
	\includegraphics[width=1\textwidth]{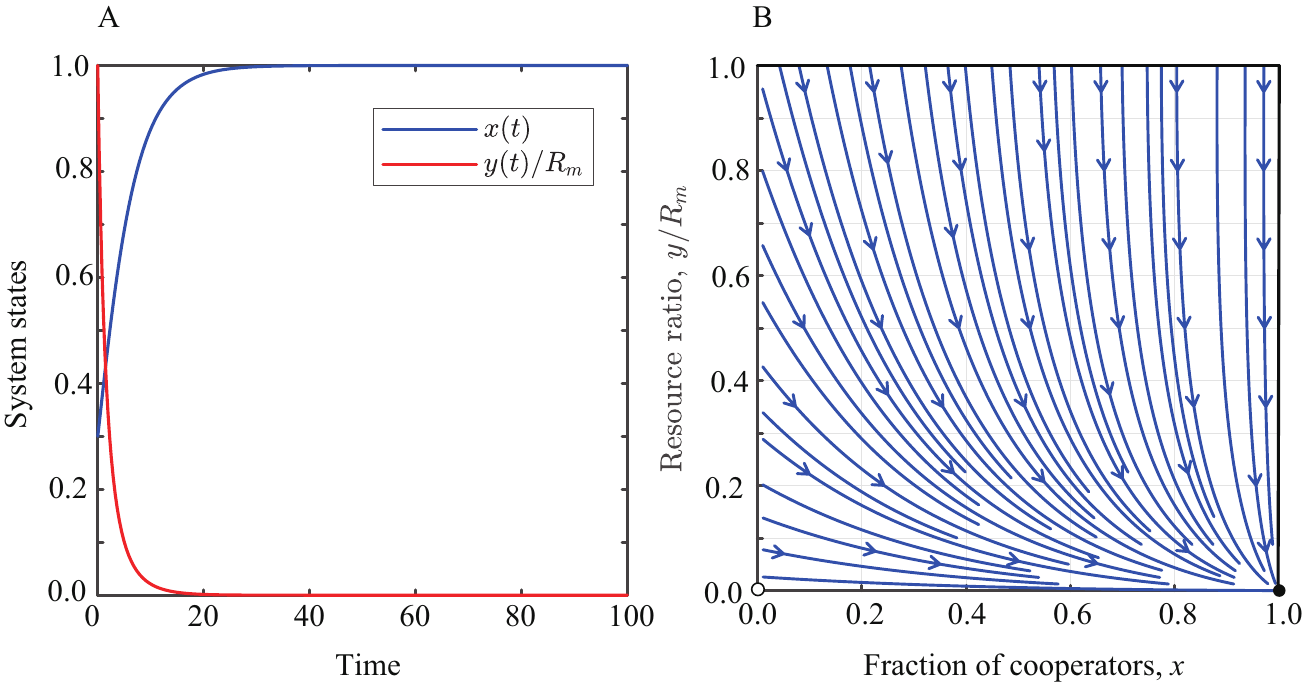}
	\caption{\textbf{The co-evolutionary dynamics of dynamic feedback games when the growth rate of the resource pool is slow}. Panel A presents the evolution of the system's state over time. Panel B shows the Phase diagram of $x-y/R_{m}$. Parameters are: $r=0.25$, $\delta=0.2$, $N=1000$, $R_m=1000$, $\alpha=0.5$, and $b_m=0.5$.}
	\label{slowfigA}
\end{figure}

\leftline{\textbf{Slowly Growing Resource Pool}}
\quad \ In the case of slow resource growth, when $0 < r < E_C < E_D$. Since $x \in [0,1]$ and $y \in [0,R_m]$, the system has exactly two boundary fixed points, $(0,0)$ and $(1,0)$. It is evident that, by analyzing the sign of the eigenvalues, $(0,0)$ is unstable, while $(1,0)$ is stable. Figure \ref{slowfigA} elegantly illustrates the dynamic interplay between the cooperation rate and the abundance of the shared resource pool over time, under the conditions where  $0 < r < E_C < E_D$. This captivating depiction is further complemented by an accompanying phase diagram, which effectively encapsulates the interaction of these variables within the specified parameter space. From the figure, an important result can be observed: when the natural growth rate $r$ is low, the reward mechanism drives participants towards full cooperation. However, even with significant rewards for cooperation, the system's resources will eventually be completely depleted, indicating that the resources are not sustainable. In this case, cooperative behavior may lose its value due to resource exhaustion. Therefore, it is advisable to avoid exploiting this type of common resource pool.\vspace*{-7pt}

\vbox{}
\leftline{\textbf{Moderately growing resource pool}}

\begin{figure}[!h]
	\centering  
	\includegraphics[width=1.0\textwidth]{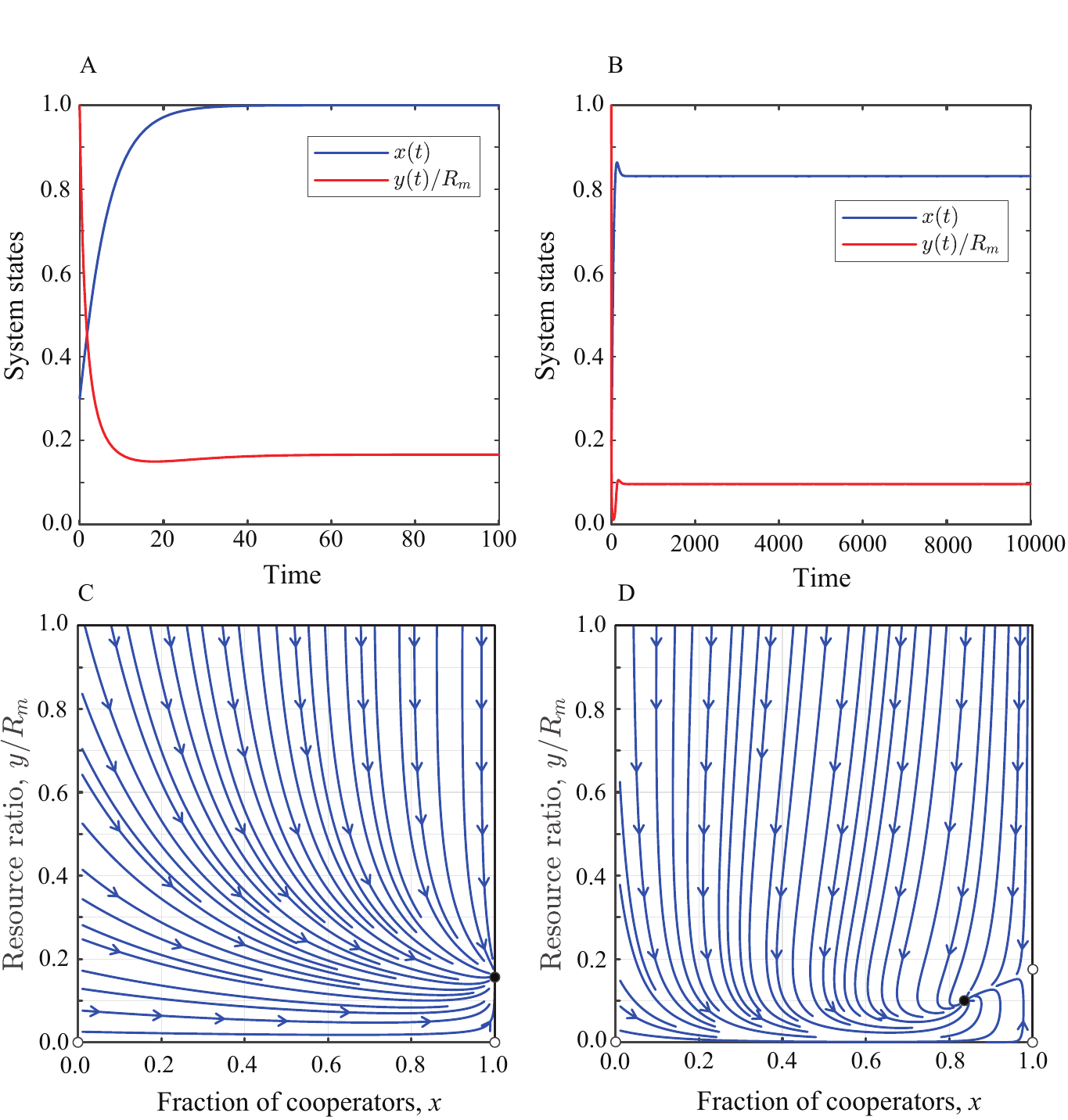}
	\caption{\textbf{The co-evolutionary dynamics of dynamic feedback games when the growth rate of the resource pool is moderate}. Panel A and C presents the evolution of the system's state over time. Panel B and D shows the Phase diagram of $x-y/R_{m}$.  The parameters for panel A and C are: $r=0.6$, $\delta=0.2$, $N=1000$, $R_m=1000$, $\alpha=0.5$, and $b_m=0.5$.  The parameters for panel B and D are: $r=0.6$, $\delta=0.02$, $N=1000$, $R_m=1000$, $\alpha=0.5$, and $b_m=0.5$.}
	\label{modfigA}
\end{figure}

\quad \ In the case of moderate resource growth, when $0 < E_C < r < E_D$. We have $N\delta - \alpha b_m + \frac{\alpha N b_m^2}{r R_m} + \frac{\alpha^2 N b_m^2}{r R_m} > 0$. We find that depending on the value of tax revenue $\delta$, the system may stabilize at different fixed points. When $\delta - \alpha b_m + \frac{\alpha N b_m^2}{r R_m} > 0$, the system has only three boundary fixed points, $\left(0,0\right)$, $\left(1,0\right)$, and $\left(1,R_m-\frac{Nb_m}{r}\right)$. By analyzing the signs of the eigenvalues, we determine that in this case, the system has a unique stable fixed point at $\left(1,R_m-\frac{Nb_m}{r}\right)$. It can be observed that when tax revenue $\delta$ is high, the system eventually leads to all participants transforming into cooperators who share the common resource pool, allowing for sustainable development of the resource pool. This intriguing scenario unfolds vividly in the premier column of Figure \ref{modfigA}.

In the event that the expression $\delta - \alpha b_m + \frac{\alpha N b_m^2}{r R_m}$ is less than zero, the system elegantly unfolds to reveal three boundary fixed points of notable significance: $\left(0, 0\right)$, $\left(1,0\right)$, and $\left(1,R_m-\frac{Nb_m}{r}\right)$, as well as one interior fixed point $\left(x^\ast, y^\ast\right)$. According to the sign analysis of eigenvalues, the first three boundary fixed points are unstable, while the last interior fixed point is stable (the stability proof for interior points can be found in Supplementary Information). It can be observed that when the tax revenue $\delta$ is low, the system will eventually reach a state of coexistence between cooperators and defectors. In this case, the resource environment can also be sustainably utilized. This situation is shown in the second column of Figure \ref{modfigA}.\vspace*{-7pt}

\vbox{}
\leftline{\textbf{Rapidly growing resource pool}}
\quad \ In the case of rapid resource growth, when $r > E_D$. We find that the choice of tax revenue $\delta$ is crucial. Depending on the different levels of tax revenue, the system may stabilize at different fixed points. However, due to the high natural growth rate of resources, the resource level in the system will always be maintained at a relatively high level for sustainable development. The stability results of the system corresponding to different levels of tax revenue are provided below.

\begin{figure}[t]
	\centering
	\includegraphics[width=1\textwidth]{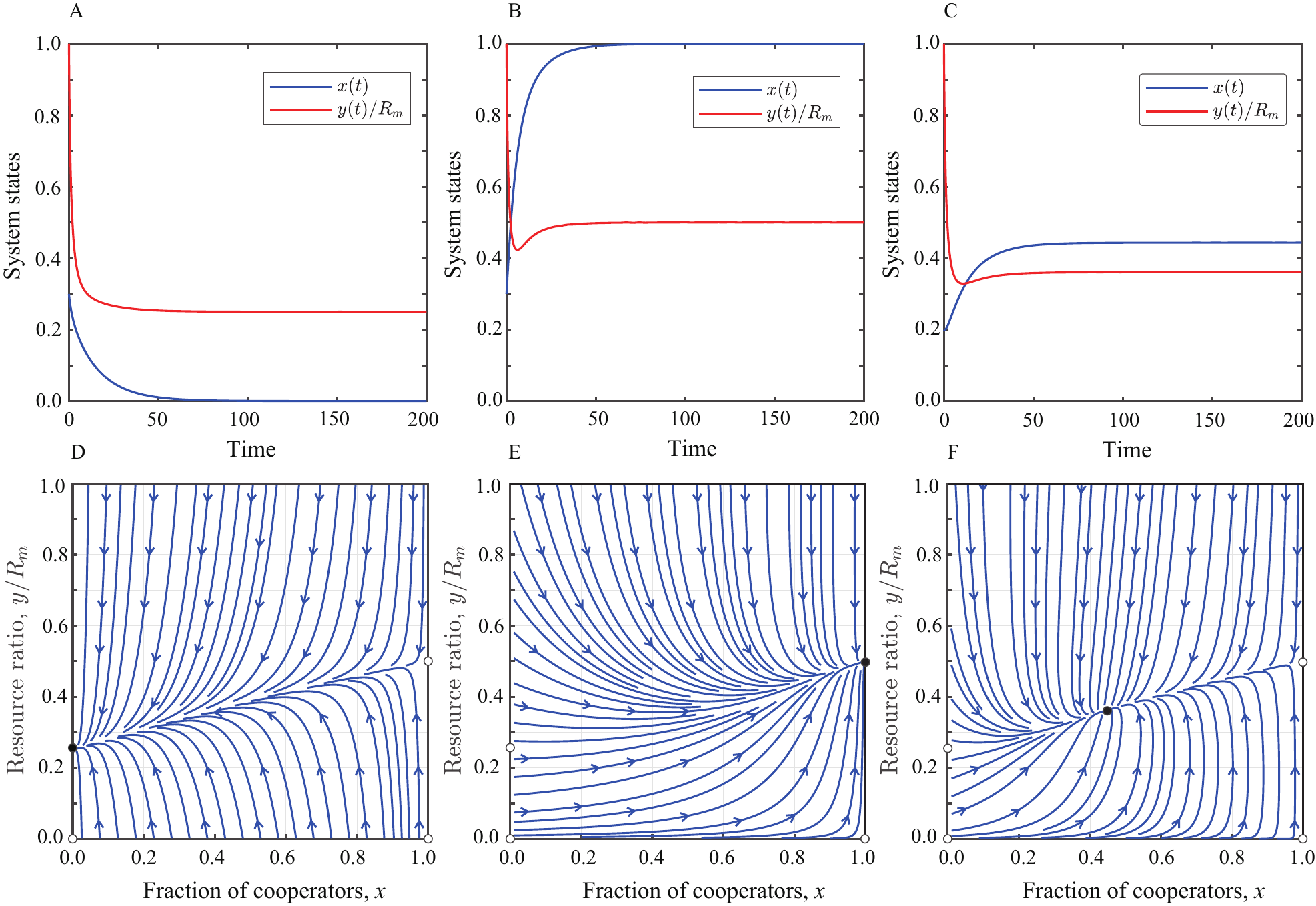}
	\caption{\textbf{The co-evolutionary dynamics of dynamic feedback games when the growth rate of the resource pool is high.} Panel A, B and C presents the evolution of the system's state over time. Panel D, E and F shows the Phase diagram of $x-y/R_{m}$.  The parameters for panel A and D are $r=1.0$, $\delta=0.00002$, $N=1000$, $R_m=1000$, $\alpha=0.5$, and $b_m=0.5$.  The parameters for panel B and E are $r=1.0$, $\delta=0.2$, $N=1000$, $R_m=1000$, $\alpha=0.5$, and $b_m=0.5$.  The parameters for panel C and F are $r=1.0$, $\delta=0.04$, $N=1000$, $R_m=1000$, $\alpha=0.5$, and $b_m=0.5$.}
	\label{repfigA}
\end{figure}

When $N\delta-\alpha b_m+\frac{\alpha N b_m^2}{r R_m}+\frac{\alpha^2 N b_m^2}{r R_m}<0$, the system elegantly unveils itself to possess merely four boundary fixed points. These points, each a testament to the system's inherent intricacy, are as follows: $\left(x,y\right)=\left(0,0\right)$, $\left(1,0\right)$, $\left(1,R_m-\frac{Nb_m}{r}\right)$, and $\left(0,R_m-\frac{Nb_m}{r}\left(1+\alpha\right)\right)$. Through our theoretical analysis, we ascertain that the fixed point $\left(0,R_m-\frac{Nb_m}{r}\left(1+\alpha\right)\right)$ is stable, while the other fixed points are unstable. This is an interesting result, as with low tax revenue $\delta$, all participants in the system eventually become defectors. However, due to the high natural growth rate of resources, the resource level in the system is maintained at a relatively high level for sustainable development. This situation is shown in the first column of Figure \ref{repfigA}.

When $\delta-\alpha b_m + \frac{\alpha N b_m^2}{r R_m}>0$, the system, in all its complexity, elegantly presents a quartet of boundary fixed points. These points, each a testament to the system's elaborate dynamics, are distinctly located at: $\left(x,y\right)=\left(0,0\right)$, $\left(1,0\right)$, $\left(1,R_m-\frac{Nb_m}{r}\right)$, $\left(0,R_m-\frac{Nb_m}{r}\left(1+\alpha\right)\right)$. It can be proven that the fixed point $\left(1,R_m-\frac{Nb_m}{r}\right)$ is stable, while the other fixed points are unstable. In this case, we can observe that when the tax revenue $\delta$ is relatively high, all participants in the system will eventually become cooperators, and the resource level in the system will be maintained at a relatively high level for sustainable development. This situation is shown in the second column of Figure \ref{repfigA}.

When $\delta-\alpha b_m + \frac{\alpha N b_m^2}{r R_m}<0<N\delta-\alpha b_m+\frac{\alpha N b_m^2}{r R_m}+\frac{\alpha^2 N b_m^2}{r R_m}$, the system has an interior fixed point $\left(x^\ast,y^\ast\right)$. Consequently, the system elegantly unfurls to reveal four boundary fixed points: $\left(x,y\right)=\left(0,0\right)$, $\left(1,0\right)$, $\left(1,R_m-\frac{Nb_m}{r}\right)$, $\left(0,R_m-\frac{Nb_m}{r}\left(1+\alpha\right)\right)$, and one interior fixed point $\left(x^\ast,y^\ast\right)$. Through analysis, we are able to ascertain that the interior fixed point $\left(x^\ast,y^\ast\right)$ is stable (for the stability of interior fixed points, please refer to Supplementary Information), while the four boundary fixed points are unstable. In this case, we can observe that the system will eventually reach a state of coexistence between cooperators and defectors. In this state, the resource level in the system will be maintained at a relatively high level for sustainable development. This situation is shown in the third column of Figure \ref{repfigA}.\vspace*{-7pt}

\vbox{}
\leftline{\textbf{Tax-punishment}}

\quad \ Drawing parallels to the scenario involving rewards, an examination of the tax-punishment strategy through the lens of equation \eqref{eq_2} unveils a similar pattern. The system, in its intricate dance of numbers, reveals four boundary fixed points: $\left(x,y\right)=\left(0,0\right)$, $\left(1,0\right)$, $\left(1,R_m-\frac{Nb_m}{r}\right)$, $\left(0,R_m-\frac{Nb_m}{r}\left(1+\alpha\right)\right)$. The detailed existence conditions and stability analysis of the interior equilibrium points are presented in the Supplementary Information. Based on the different natural growth rate $r$, we will discuss the following three cases.

\begin{figure}[t]
	\centering
	\includegraphics[width=1\textwidth]{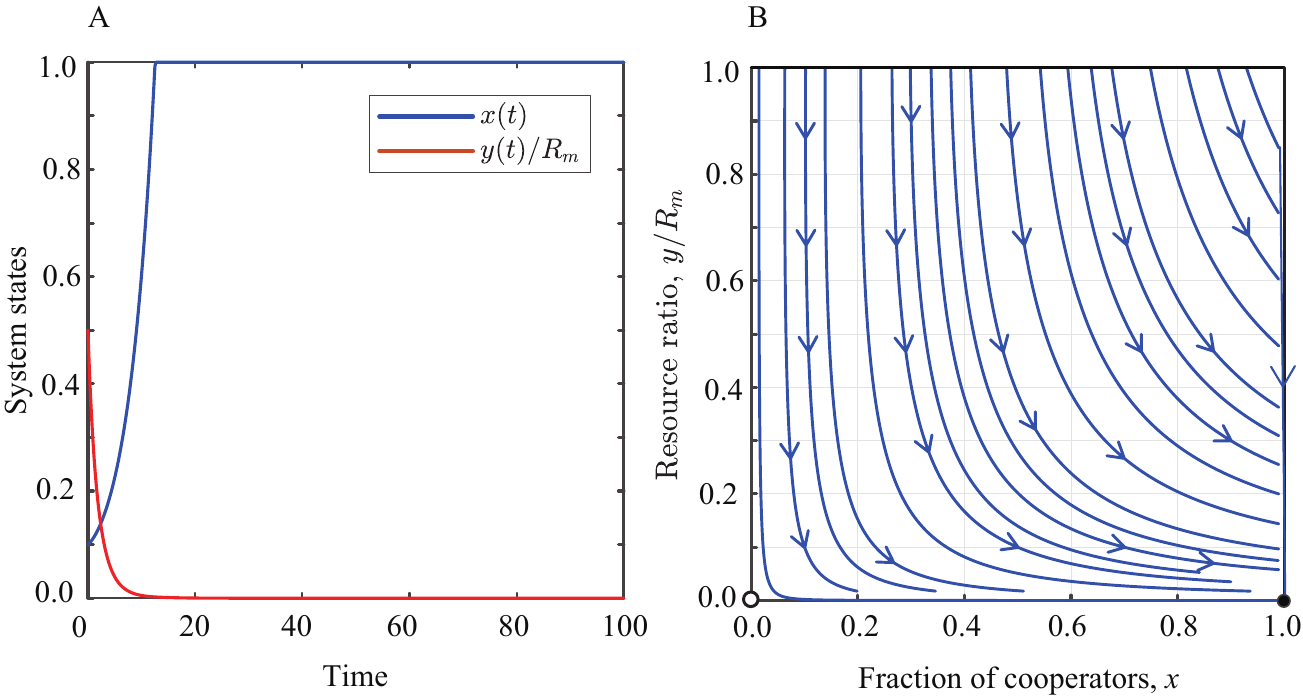}
	\caption{\textbf{The co-evolutionary dynamics of coupled systems, under sscenarios of a slowly growing resource pool, incorporate tax-punishment.} Panel A presents the evolution of the system's states over time. Panel B shows the Phase diagram of $x-y/R_{m}$. Parameters are $r=0.25$, $\delta=0.2$, $N=1000$, $R_m=1000$, $\alpha=0.5$, and $b_m=0.5$.}
	\label{fig4}
\end{figure}

\leftline{\textbf{Slowly Growing Resource Pool}}
\quad \ In the case of slow resource growth, that is $0 < r < E_C < E_D$, the system has two and only two equilibrium points: $(0, 0)$ and $(1, 0)$. Obviously, by analyzing the signs of the eigenvalues, we can determine that $(0, 0)$ is unstable, while $(1, 0)$ is stable. Figure \ref{fig4} elegantly unfurls the tapestry of time, charting the evolutionary trajectory of the proportion of cooperators and the abundance of the common resource pool. This intricate dance of variables plays out under the condition where $0 < r < E_C < E_D$. The figure also gracefully presents the corresponding phase diagram, a testament to the system's dynamism under the specified parameters. This visual representation serves as a window into the system's complex interplay of factors, illuminating the intricate relationships that govern its behavior. Similar to the case of tex-reward, when the natural growth rate $r$ is low, the punishing mechanism drives the participants towards complete cooperation. However, even with significant penalties imposed on defectors, the resources will eventually be completely depleted. In this scenario, high levels of cooperation become meaningless.\vspace*{-7pt}

\begin{figure}[t]
	\centering
	\includegraphics[width=1\textwidth]{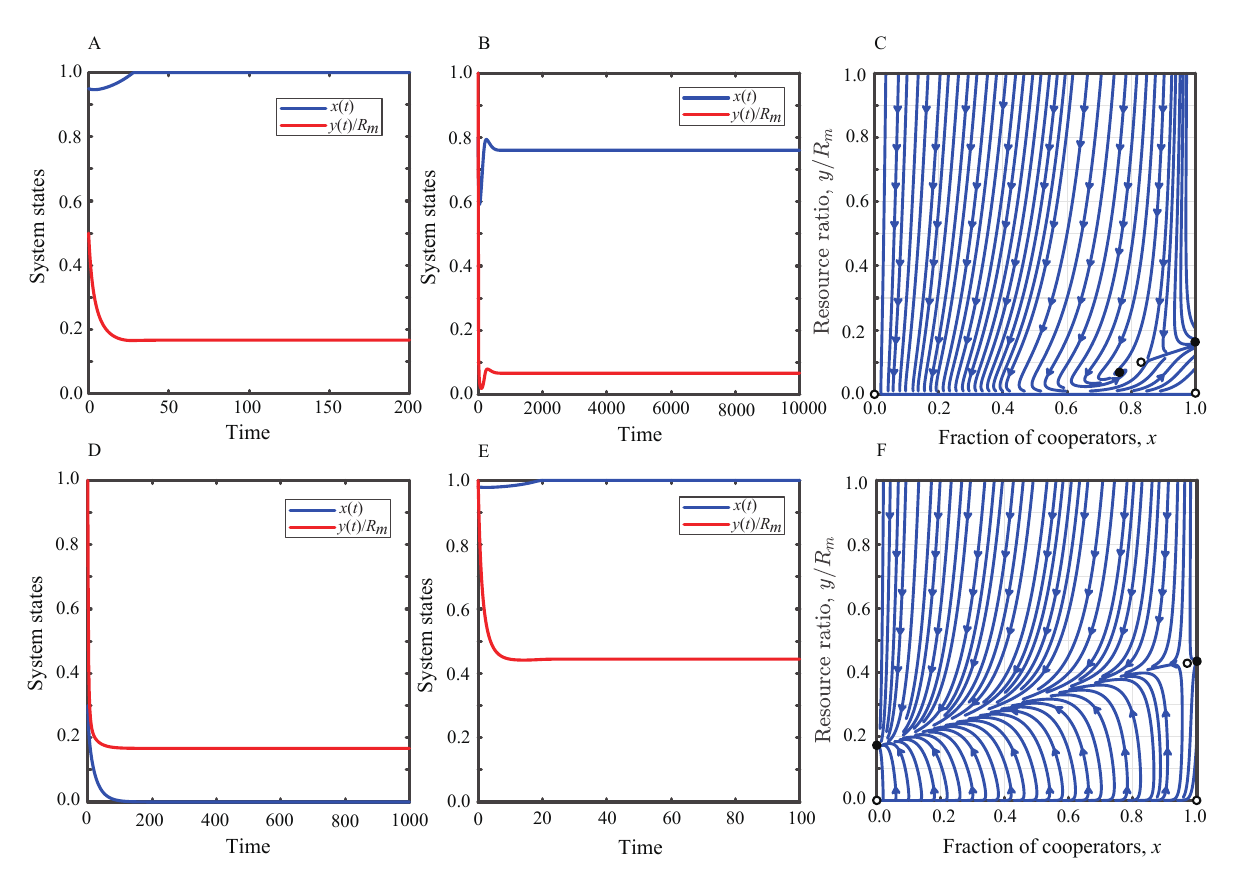}
	\caption{\textbf{The coupled system can produce bistable outcomes.} Panels A, B, D and E presents the evolution of the system's state over time. Panel C and F shows the Phase diagram of $x-y/R_{m}$.  The parameters for panels A, B, and C are $r=0.6$, $\delta=0.004$, $N=1000$, $R_m=1000$, $\alpha=0.5$, and $b_m=0.5$. The parameters for panels D, E and F are $r=0.9$, $\delta=0.003$, $N=1000$, $R_m=1000$, $\alpha=0.5$, and $b_m=0.5$.}
	\label{modfigpA}
\end{figure}

\vbox{}
\leftline{\textbf{Moderately growing resource pool}}
\quad \ When resource growth is moderate, specifically $0 < E_C < r < E_D$. For the boundary fixed point, it is evident that $\left(1,R_m-\frac{Nb_m}{r}\right)$ is stable when $-N\delta+\alpha b_m-\frac{\alpha N b_m^2}{r R_m}<0$ (More detailed theoretical analysis can be found in the Supplementary Information). Now we provide a numerical example. We find that the system has three boundary fixed points $\left(0, 0\right)$, $\left(1, 0\right)$, $\left(1,R_m-\frac{Nb_m}{r}\right)$ and two interior fixed points $\left(x_1^\ast,y_1^\ast\right)$, $\left(x_2^\ast,y_2^\ast\right)$. Let's assume $x_1^\ast<x_2^\ast$. In this situation, the system has two stable fixed points $\left(1,R_m-\frac{Nb_m}{r}\right)$ and $\left(x_1^\ast,y_1^\ast\right)$, whereas the remaining fixed points are all unstable. This scenario is illustrated in the first row of Figure \ref{modfigpA}. We find that there are five equilibrium points in the phase plane, among which there are two stable equilibria. One is located inside the phase plane, indicating that cooperation can be maintained at a relatively high level and the resource level can also remain constant. The other is located on the right boundary of the phase plane, implying that all individuals in the group choose cooperative behavior while resources can be sustained. It should be noted that the emergence of one of the evolutionary outcomes depends on the initial levels of cooperation and resources.

\vbox{}
\leftline{\textbf{Rapidly growing resource pool}}

\quad \ When the natural growth rate is high, specifically $r>E_D$. Clearly, when $-N\delta+\alpha b_m-\frac{\alpha N b_m^2}{r R_m}<0$, the fixed point $\left(1, R_m-\frac{Nb_m}{r}\right)$ is stable; when $\delta-\alpha b_m+\frac{\alpha N b_m^2}{r R_m}+\frac{\alpha^2 N b_m^2}{r R_m}<0$, the fixed point $\left(0, R_m-\frac{Nb_m\left(1+\alpha\right)}{r}\right)$ is stable (Detailed theoretical analysis is presented in the Supplementary Information). Here, we provide a numerical example that ensures the existence of an unstable interior equilibrium point in the coupled system. As depicted in the second row of Figure \ref{modfigpA}, there exist five equilibrium points in the phase plane: $\left(0,0\right)$, $\left(1,0\right)$, $\left(1,R_m-\frac{Nb_m}{r}\right)$, $\left(0,R_m-\frac{Nb_m\left(1+\alpha\right)}{r}\right)$, and one interior fixed point $\left(x^\ast,y^\ast\right)$. In this situation, the two boundary fixed points $\left(1,R_m-\frac{Nb_m}{r}\right)$ and $\left(0,R_m-\frac{Nb_m\left(1+\alpha\right)}{r}\right)$ are stable, but the interior fixed point $\left(x^\ast,y^\ast\right)$ is unstable. Similar to the scenario with a moderate natural growth rate, the system displays bistable outcomes, underscoring the importance of the selection of initial conditions, as the system's stability hinges on these chosen values. Despite the potential for the system to ultimately evolve towards either complete cooperation or complete defection, the high natural growth rate of the common resource pool enables the system to consistently maintain a non-zero level of resources.

\begin{figure}[t]
	\centering
	\includegraphics[width=1\textwidth]{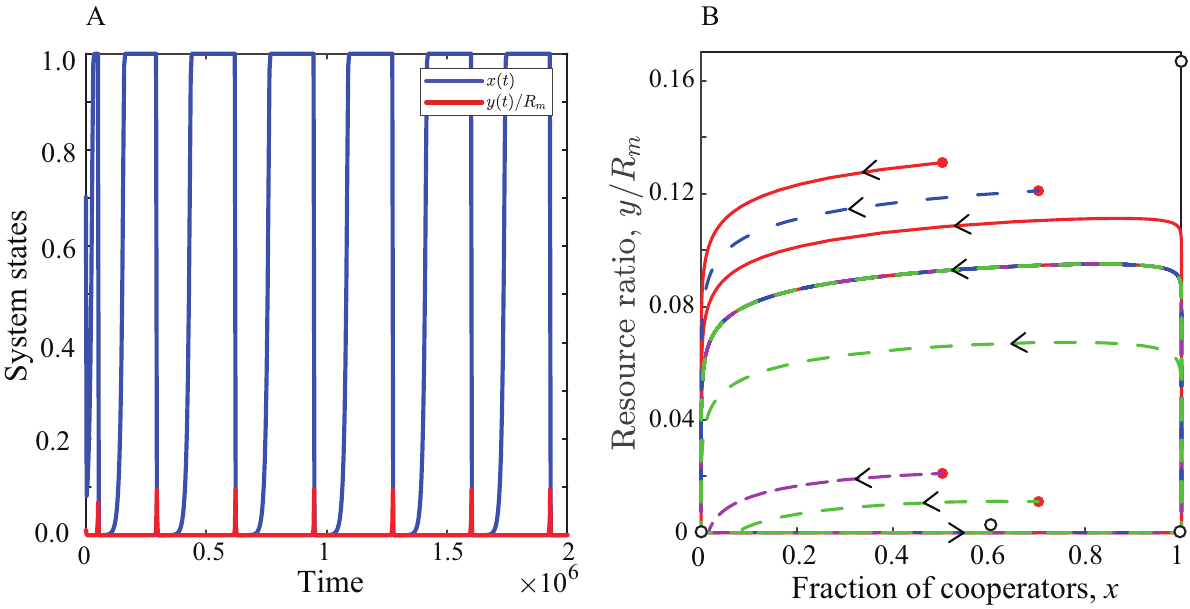}
	\caption{\textbf{The coupled system can generate evolutionary oscillatory dynamics.} Panel A presents the evolution of the system's state over time. Panel B shows the Phase diagram of $x-y/R_{m}$. Parameters are $r=0.006$, $\delta=0.0001$, $N=10$, $R_m=1000$, $\alpha=0.5$, and $b_m=0.5$.}
	\label{bwd}
\end{figure}

It is worth noting that the coupled system can generate limit cycle dynamics. Figure \ref{bwd} shows the representative evolution of cooperation level and abundance of the common resource pool over time, as well as the phase diagram of $x-y/Rm$ under the condition of $E_C< r < E_D$ (The theoretical analysis can be found in Case 1 of the Supplementary Information). We observe four equilibrium points in the phase plane, of which three are boundary equilibrium points and one is an interior equilibrium point, all of which are unstable. Trajectories originating from any initial point in the phase plane ultimately converge to a limit cycle within the plane (see Figure \ref{bwd}B). This implies that the frequency of cooperators and the level of resources exhibit oscillating evolutionary dynamics. The trajectory of system states over time in Figure \ref{bwd}A also verifies this result.

\leftline{\textbf{Monte Carlo Simulations}}
\quad \ In our quest to corroborate the efficacy of our model, we embarked on a series of Monte Carlo numerical simulations. The application of Monte Carlo methodologies within the framework of system dynamics serves as a powerful tool to augment our comprehension of system behavior, and more importantly, to bolster our predictive prowess. This, in turn, equips decision-makers with the capacity to devise more potent strategies and to appraise the risks and uncertainties intertwined with their decisions with greater precision. In the forthcoming section, we will elucidate the robustness and effectiveness of our dual-model framework, as substantiated by our comprehensive Monte Carlo simulations.

\leftline{\textbf{Tax-reward}}
\quad \ For each individual $i$ in the population, they are required to pay taxes $\delta$, which will be evenly distributed as reward to cooperators. Each individual has only two choices: to become a cooperator (denoted as $s_i=1$) or become a defector (denoted as $s_i=0$). If individual $i$ is a cooperator, his/her payoff is $F_i=b_my(t)/R_m-\delta+N\delta/\sum_{i=1}^N s_i$. Otherwise, if individual $i$ is a defector, his/her payoff is $F_i=\left(1+\alpha\right) b_my(t)/R_m-\delta$. The update of the common resource level $y(t)$ is governed by the following rule:
\begin{eqnarray}
	y(t+1)=y(t)+r y(t)\left[1-\frac{y(t)}{R_{m}}\right]-\sum_{i=1}^N \left[s_{i} \frac{y(t) b_{m}}{R_{m}}+\left(1-s_{i}\right) \frac{y(t) b_{m}(1+\alpha)}{R_{m}}\right].
\end{eqnarray}

The selection of individual strategies involves each individual $i$ learning from a random individual $j$. If the payoff $F_i$ for individual $i$ under his/her chosen strategy is lower than the payoff $F_j$ for individual $j$, then individual $i$ will adopt the strategy of $j$ with a probability $q$. Otherwise, individual $i$ will remain with their current strategy. We refer to the probability $q$ as the transition probability.
\begin{eqnarray}
	q=\frac{F_{i}-F_{j}}{M},
\end{eqnarray}
where $M$ quantifies the uncertainty in strategy adoption. 

\begin{figure}[t]
	\centering  
	\includegraphics[width=1\textwidth]{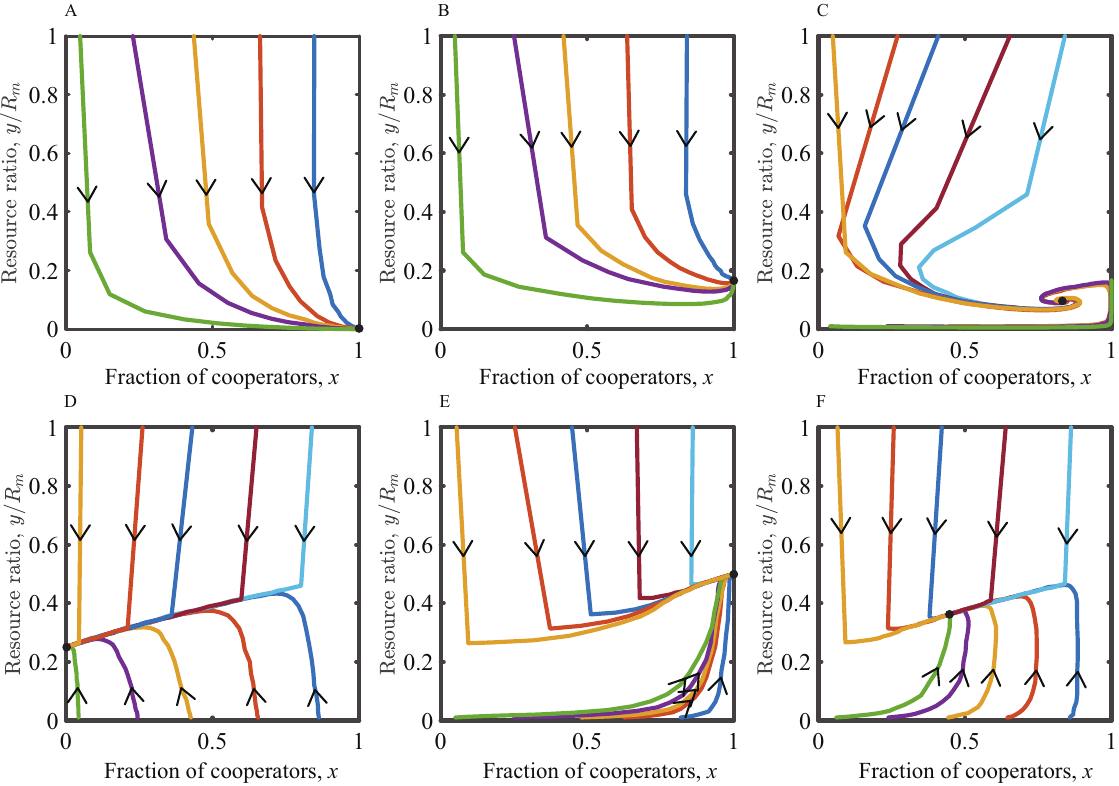}
	\caption{\textbf{The results of Monte Carlo simulations for the coupled system with tax-reward under three distinct scenarios.} The solid black dots in the figures represent the calculated stable fixed points. Panel A displays the simulation outcomes when the resource growth rate is low, while Panels B and C present two different sets of results for scenarios with a moderate resource growth rate. Panels D, E, and F exhibit three distinctly varied simulation outcomes under conditions of a high resource growth rate. The parameters for panel A are $M=1$, $r=0.25$, $\delta=0.2$, $N=1000$, $R_m=1000$, $\alpha=0.5$, and $b_m=0.5$. The parameters for panel B are $M=1$, $r=0.6$, $\delta=0.2$, $N=1000$, $R_m=1000$, $\alpha=0.5$, and $b_m=0.5$. The parameters for panel C are $M=0.1$, $r=0.6$, $\delta=0.02$, $N=1000$, $R_m=1000$, $\alpha=0.5$, and $b_m=0.5$. The parameters for panel D are $M=1$, $r=1$, $\delta=0.00002$, $N=1000$, $R_m=1000$, $\alpha=0.5$, and $b_m=0.5$. The parameters for panel E are $M=1$, $r=1$, $\delta=0.2$, $N=1000$, $R_m=1000$, $\alpha=0.5$, and $b_m=0.5$. The parameters for panel F are $M=1$, $r=1$, $\delta=0.04$, $N=1000$, $R_m=1000$, $\alpha=0.5$, and $b_m=0.5$.}
	\label{mtkl}
\end{figure}

Figure \ref{mtkl} illustrates the results of the evolution of a resource pool over time under different scenarios: slow growth, moderate growth, and rapid growth, as obtained from Monte Carlo simulations. From the figures, it can be observed that the system eventually stabilizes at the theoretically calculated fixed points. In the case of slow resource growth, the magnitude of rewards is insufficient to sustain resource development. Although all participants eventually transition to cooperators, the resources will be completely depleted (see Figure \ref{mtkl} A). Under moderate resource growth, consistent results with the theoretical expectations are obtained, where the stability of the system is determined by the magnitude of rewards (see Figure \ref{mtkl} B and C). In the scenario of rapid resource growth, depending on the reward levels, the system may reach a state of complete defection (see Figure \ref{mtkl} D), complete cooperation (see Figure \ref{mtkl} E), or coexistence of cooperation and defection (see Figure \ref{mtkl} F). Irrespective of which among the aforementioned three scenarios manifests, the resources in the public pool consistently maintain a stable level.

\vbox{}
\leftline{\textbf{Tax-punishment}}
\quad \ For the tax-punishment strategy, if individual $i$ is a cooperator, his/her payoff is $F_i=b_my(t)/R_m-\delta$. Otherwise, if individual $i$ is a defector, his/her payoff is $F_i=\left(1+\alpha\right)b_my(t)/R_m-\delta-N\delta/\sum_{i=1}^{N}(1-s_i)$. We still adopt the same updating strategy for the common resource pool and individual strategies as in the case of reward.

\begin{figure}[t]
	\centering  
	\includegraphics[width=1\textwidth]{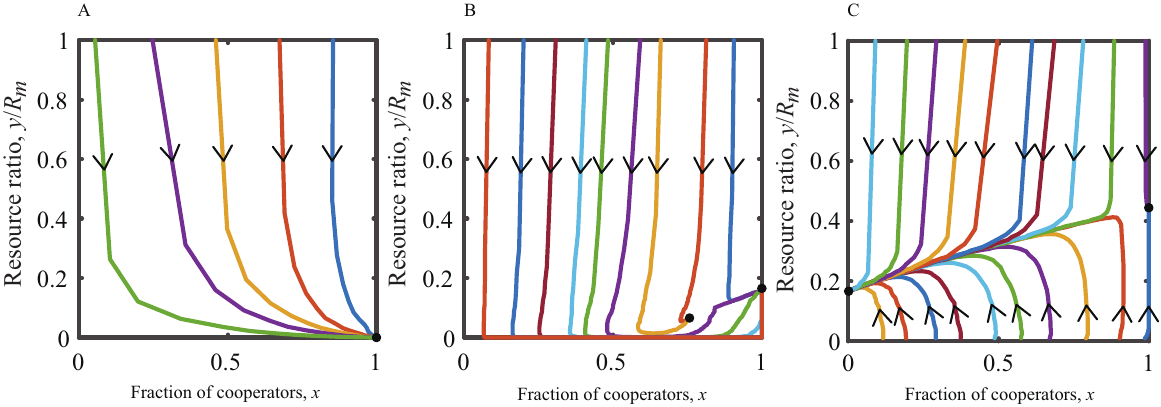}
	\caption{\textbf{Monte Carlo simulation results of the co-evolutionary game model under the tax penalty scenario when considering low, moderate, and high resource growth rates, respectively.} The solid black dots in the figures represent the calculated stable fixed points. The parameters for panel A are $M=1$, $r=0.25$, $\delta=0.2$, $N=1000$, $R_m=1000$, $\alpha=0.5$, and $b_m=0.5$. The parameters for panel B are $M=2$, $r=0.6$, $\delta=0.004$, $N=1000$, $R_m=1000$, $\alpha=0.5$, and $b_m=0.5$. The parameters for panel C are $M=1$, $r=0.9$, $\delta=0.003$, $N=1000$, $R_m=1000$, $\alpha=0.5$, and $b_m=0.5$. }
	\label{mtkl1}
\end{figure}

In Figure \ref{mtkl1}, we present the agent-based simulation results considering scenarios with low, moderate, and high growth rates of the public resource pool, respectively, when taking into account tax-punishment. Our simulation results corroborate the theoretical findings: when the growth rate of the public pool is low, despite all individuals in the group ultimately opting for a cooperative strategy, the resources inevitably reach a state of complete depletion (see Figure \ref{mtkl1}A). This implies that even universal cooperation cannot overcome the tragedy of the commons induced by an excessively low resource growth rate. When the growth rate of the public pool is moderate, depending on the initial conditions, the system either stabilizes in a state where cooperators and defectors coexist sustainably with resources, or it stabilizes in a state where all individuals choose cooperation accompanied by sustainable resources (see Figure \ref{mtkl1}B). Either outcome results in a situation where both cooperation and resources are maintained. When the growth rate of the public pool is high, another bistable outcome can emerge: regardless of whether the population ultimately evolves into a state of complete cooperation or complete defection, resources can be sustained (see Figure \ref{mtkl1}C). Due to the high natural growth rate of resources in this case, the resource level of the system can be sustained at a higher level and be developed sustainably.
\vspace*{-5pt}

\vbox{}
\leftline{\textbf{Conclusions}}
\quad \ 

The pivotal importance of sustainable resource development is gaining heightened acknowledgment in our contemporary society. Consequently, forecasting the enduring impacts of individual actions on resource stewardship has ascended to the forefront as a paramount challenge. Existing research has explored the complex relationship between resource environments and stakeholder behaviors, leading to the formulation of sophisticated models \cite{traulsen2023future,yc1,first2,second1,second2,third1,wang2020eco}. These models integrate diverse factors such as incentives \cite{chenpcb2018}, regulatory norms \cite{ostrom1990governing}, and sociocultural influences \cite{leach1999environmental}, providing a comprehensive perspective on resource utilization dynamics. Although relevant research has been initiated, the understanding of how to better design incentives to regulate individual behaviors for promoting resource sustainability remains in a nascent stage.

We have proposed new incentive strategies and established tax-reward and tax-punishment models to more accurately describe the reality of resource development. We have found that the effects generated by rewarding cooperators and punishing defectors are not uniform. Specifically, in scenarios characterized by a low natural growth rate of resources, neither the tax-reward nor the tax-punishment model is capable of achieving sustainable resource development solely through rewarding cooperators or punishing defectors. Regardless of the magnitude of rewards given to cooperators or punishments imposed on defectors, although participants in the system tend to become all cooperators, the resources in the environment will eventually deplete. Therefore, for resources with low natural growth rates, more effective protection measures should be implemented, such as avoiding exploitation, to prevent eventual resource depletion. For resources with high natural growth rates, our research results in the tax-reward model have demonstrated that sustainable resource development can be ensured. Here the magnitude of rewards does not affect the sustainability of resources but only the total resource level in the steady state. This provides important insights for practical applications: by adjusting the size of the tax payment $\delta$, the system can eventually stabilize at a higher resource level. Regarding the tax-punishment model, our findings suggest the potential for bistability within the system. The stability of the fixed points is not solely tethered to the magnitude of the tax payment, denoted as $\delta$, but also intricately intertwined with the choice of the initial state. This dual dependency adds a layer of complexity to the system's dynamics. In this case, choosing an appropriate initial state is also critical in maintaining a high resource level. For a resource pool with moderate natural growth rates, our research results in the tax-reward model show that the system only exhibits one stable fixed point. In reality, we need to adjust the tax revenue $\delta$ to maintain a higher resource level in the system. In our exploration of the tax-punishment model, we have discovered that the system may display bistability, suggesting that the system's stability is contingent upon the strength of punishment and the selection of the initial states. However, it's important to acknowledge that in the tax-punishment model, scenarios may arise where stable fixed points are absent. For instance, this paper presents an example of a limit cycle, wherein the system fails to achieve a stable state and instead undergoes periodic cycling. This phenomenon could stem from intricate internal interactions and feedback loops within the system. This insight is crucial for gaining a more profound understanding of the system's dynamic behavior and for predicting future evolutionary trends.

Our research results provide important theoretical foundations for decision-making in practical contexts, aiding us in better resource management and utilization to achieve sustainable development. However, we also recognize that the model has its inherent limitations. For instance, we assume that all individuals are perfectly rational, while in reality, individual decision-making may be influenced by various factors. Therefore, future research can further expand and optimize our models, such as introducing more real-life factors, considering irrationality in individual behavior \cite{hua2023evolution,sarkar2023managing,hilbe2018evolution}, accounting for resource growth rates varying with the current resource level $y$, or incorporating the migration of individuals \cite{wang2016effects,shu2023determinants}. It would also be beneficial to compare the actual cost-effectiveness of the tax-reward and tax-punishment mechanisms proposed in our current work. Optimal institutional incentives for promoting cooperation, as studied in recent works within the context of public goods and prisoner's dilemma games, have shown that adaptive rewarding strategies considering population statistics, such as the threshold number of cooperators in the population, or local network properties, like the number of cooperators in the network, could be more cost-efficient in promoting cooperation \cite{han2022institutional,duong2021cost,wang2019exploring,cimpeanu2023does,duong2023cost}. Exploring the applicability of such adaptive approaches in our current context, depending on the resource growth rate and resource dynamics, could provide further insights into effective incentive mechanisms for sustainable resource management.

\clearpage

\noindent \textbf{Ethics.} \\
This work did not require ethical approval from a human subject or animal welfare committee.

\noindent \textbf{Acknowledgments} \\
This research was supported by the National Natural Science Foundation of China (Grant No. 62306243), the Natural Science Foundation of Shaanxi (Grant Nos. 2023-JC-QN-0791, 2023-JC-YB-083, 23JK0693, and 23JK0692) and the Fundamental Research Funds of the Central Universities of China (Grants No. 2452022012, 2452022144).

\noindent \textbf{Data accessibility}\\
All pertinent analysis has been included in the Supplementary Information. Source code is available at the Dryad, Dataset, https://doi.org/10.5061/dryad.pnvx0k6w6.

\noindent \textbf{Declaration of AI use}\\
We have not used AI-assisted technologies in
creating this article.

\noindent \\ \textbf{Author contributions} \\
L. W.: investigation, methodology, writing—original draft, writing—review and editing; Y. L.: investigation, methodology; R. G.: investigation; L. Z.: writing—original draft; L. L: conceptualization, formal analysis, investigation, supervision, methodology, writing—original draft, writing—review and editing; S. H.: formal analysis, methodology, supervision, writing—original draft, writing—review and editing.

\noindent \\ \textbf{Competing financial interests} \\
We declare we have no competing interests.

\newpage

\begin{thebibliography}{54}

\bibitem{plus1}
Cox M, Arnold G, Tom{\'a}s SV. 2010  A review of design principles for
  community-based natural resource management. {\em Ecology and Society}
  \textbf{15}.

\bibitem{plus2}
Hauser OP, Rand DG, Peysakhovich A, Nowak MA. 2014  Cooperating with the
  future. {\em Nature} \textbf{511}, 220--223.

\bibitem{plus3}
Sugiarto HS, Lansing JS, Chung NN, Lai C, Cheong SA, Chew LY. 2017  Social
  cooperation and disharmony in communities mediated through common pool
  resource exploitation. {\em Physical Review Letters} \textbf{118}, 208301.

\bibitem{gg1}
Ringsmuth AK, Lade SJ, Schl{\"u}ter M. 2019  Cross-scale cooperation enables
  sustainable use of a common-pool resource. {\em Proceedings of the Royal
  Society B: Biological Sciences} \textbf{286}, 20191943.

\bibitem{liu2023coevolutionary}
Liu L, Chen X, Szolnoki A. 2023  Coevolutionary dynamics via adaptive feedback
  in collective-risk social dilemma game. {\em eLife} \textbf{12}, e82954.

\bibitem{rl2}
Bellomo N, Liao J, Quaini A, Russo L, Siettos C. 2023  Human behavioral crowds:
  Review, critical analysis, and research perspectives. {\em Mathematical
  Models and Methods in Applied Sciences} \textbf{33}, 1611--1659.

\bibitem{third2}
Tanimoto J. 2019  Evolutionary games with sociophysics. {\em Evolutionary
  Economics} \textbf{17}.

\bibitem{rll1}
Hendry AP, Gotanda KM, Svensson EI. 2017  Human influences on evolution, and
  the ecological and societal consequences. .

\bibitem{gdfm1}
Gilles JL, Jamtgaard K. 1981  Overgrazing In Pastoral Areas The Commons
  Reconsidered. {\em Sociologia Ruralis} \textbf{21}, 129--141.

\bibitem{gdfm2}
Wilson AD, MacLeod ND. 1991  Overgrazing: present or absent?. {\em Rangeland
  Ecology \& Management/Journal of Range Management Archives} \textbf{44},
  475--482.

\bibitem{ref21}
Mysterud A. 2006  The concept of overgrazing and its role in management of
  large herbivores. {\em Wildlife Biology} \textbf{12}, 129--141.

\bibitem{ref19}
Kalantari Z, Ferreira CSS, Page J, Goldenberg R, Olsson J, Destouni G. 2019
  Meeting sustainable development challenges in growing cities: Coupled
  social-ecological systems modeling of land use and water changes. {\em
  Journal of Environmental Management} \textbf{245}, 471--480.

\bibitem{ref16}
Pauly D, Christensen V, Gu{\'e}nette S, Pitcher TJ, Sumaila UR, Walters CJ,
  Watson R, Zeller D. 2002  Towards sustainability in world fisheries. {\em
  Nature} \textbf{418}, 689--695.

\bibitem{ref17}
Kraak SB. 2011  Exploring the ‘public goods game’model to overcome the
  tragedy of the commons in fisheries management. {\em Fish and Fisheries}
  \textbf{12}, 18--33.

\bibitem{wang2020communicating}
Wang Z, Jusup M, Guo H, Shi L, Ge{\v{c}}ek S, Anand M, Perc M, Bauch CT, Kurths
  J, Boccaletti S et~al.. 2020  Communicating sentiment and outlook reverses
  inaction against collective risks. {\em Proceedings of the National Academy
  of Sciences} \textbf{117}, 17650--17655.

\bibitem{hz1}
Connor RC. 2010  Cooperation beyond the dyad: on simple models and a complex
  society. {\em Philosophical Transactions of the Royal Society B: Biological
  Sciences} \textbf{365}, 2687--2697.

\bibitem{hz2}
Sasaki T, Uchida S. 2013  The evolution of cooperation by social exclusion.
  {\em Proceedings of the Royal Society B: Biological Sciences} \textbf{280},
  20122498.

\bibitem{hz3}
Brosnan SF, Salwiczek L, Bshary R. 2010  The interplay of cognition and
  cooperation. {\em Philosophical Transactions of the Royal Society B:
  Biological Sciences} \textbf{365}, 2699--2710.

\bibitem{capraro2021mathematical}
Capraro V, Perc M. 2021  Mathematical foundations of moral preferences. {\em
  Journal of the Royal Society interface} \textbf{18}, 20200880.

\bibitem{sigmund2010social}
Sigmund K, De~Silva H, Traulsen A, Hauert C. 2010  Social learning promotes
  institutions for governing the commons. {\em Nature} \textbf{466}, 861--863.

\bibitem{hilbe2010incentives}
Hilbe C, Sigmund K. 2010  Incentives and opportunism: from the carrot to the
  stick. {\em Proceedings of the Royal Society B: Biological Sciences}
  \textbf{277}, 2427--2433.

\bibitem{ref4}
Fehr E, G{\"a}chter S. 2000  Cooperation and punishment in public goods
  experiments. {\em American Economic Review} \textbf{90}, 980--994.

\bibitem{ref5}
Henrich J. 2006  Cooperation, punishment, and the evolution of human
  institutions. {\em Science} \textbf{312}, 60--61.

\bibitem{ref7}
Sigmund K, De~Silva H, Traulsen A, Hauert C. 2010  Social learning promotes
  institutions for governing the commons. {\em Nature} \textbf{466}, 861--863.

\bibitem{ref8}
Vasconcelos VV, Santos FC, Pacheco JM. 2013  A bottom-up institutional approach
  to cooperative governance of risky commons. {\em Nature Climate Change}
  \textbf{3}, 797--801.

\bibitem{ref10}
Johnson S. 2015  Escaping the tragedy of the commons through targeted
  punishment. {\em Royal Society Open Science} \textbf{2}, 150223.

\bibitem{ref12}
Szolnoki A, Perc M. 2015  Antisocial pool rewarding does not deter public
  cooperation. {\em Proceedings of the Royal Society B: Biological Sciences}
  \textbf{282}, 20151975.

\bibitem{ref13}
Perc M, Jordan JJ, Rand DG, Wang Z, Boccaletti S, Szolnoki A. 2017  Statistical
  physics of human cooperation. {\em Physics Reports} \textbf{687}, 1--51.

\bibitem{jl1}
Dong Y, Sasaki T, Zhang B. 2019  The competitive advantage of institutional
  reward. {\em Proceedings of the Royal Society B: Biological Sciences}
  \textbf{286}, 20190001.

\bibitem{wang2023emergence}
Wang Z, Song Z, Shen C, Hu S. 2023  Emergence of punishment in social dilemma
  with environmental feedback. In {\em Proceedings of the AAAI Conference on
  Artificial Intelligence} vol.~37 pp. 11708--11716.

\bibitem{jl2}
Gneezy A, Fessler DM. 2012  Conflict, sticks and carrots: war increases
  prosocial punishments and rewards. {\em Proceedings of the Royal Society B:
  Biological Sciences} \textbf{279}, 219--223.

\bibitem{chenpcb2018}
Chen X, Szolnoki A. 2018  Punishment and inspection for governing the commons
  in a feedback-evolving game. {\em PLoS Computational Biology} \textbf{14},
  e1006347.

\bibitem{ref23}
Chen X, Gross T, Dieckmann U. 2013  Shared rewarding overcomes defection traps
  in generalized volunteer's dilemmas. {\em Journal of Theoretical Biology}
  \textbf{335}, 13--21.

\bibitem{hua2024coevolutionary}
Hua S, Liu L. 2024  Coevolutionary dynamics of population and institutional
  rewards in public goods games. {\em Expert Systems with Applications}
  \textbf{237}, 121579.

\bibitem{ref24}
Chen X, Sasaki T, Br{\"a}nnstr{\"o}m {\AA}, Dieckmann U. 2015  First carrot,
  then stick: how the adaptive hybridization of incentives promotes
  cooperation. {\em Journal of the Royal Society Interface} \textbf{12},
  20140935.

\bibitem{sun2021combination}
Sun W, Liu L, Chen X, Szolnoki A, Vasconcelos VV. 2021  Combination of
  institutional incentives for cooperative governance of risky commons. {\em
  iScience} \textbf{24}, 102844.

\bibitem{ref25}
Schuster P, Sigmund K. 1983  Replicator dynamics. {\em Journal of Theoretical
  Biology} \textbf{100}, 533--538.

\bibitem{ref26}
Hauert C, De~Monte S, Hofbauer J, Sigmund K. 2002  Replicator dynamics for
  optional public good games. {\em Journal of Theoretical Biology}
  \textbf{218}, 187--194.

\bibitem{ref27}
Nowak MA, Sigmund K. 2004  Evolutionary dynamics of biological games. {\em
  Science} \textbf{303}, 793--799.

\bibitem{ref28}
Weitz JS, Eksin C, Paarporn K, Brown SP, Ratcliff WC. 2016  An oscillating
  tragedy of the commons in replicator dynamics with game-environment feedback.
  {\em Proceedings of the National Academy of Sciences} \textbf{113},
  E7518--E7525.

\bibitem{ref35}
Hofbauer J, Sigmund K. 1998 {\em Evolutionary games and population dynamics}.
Cambridge University Press.

\bibitem{ref36}
Tsoularis A, Wallace J. 2002  Analysis of logistic growth models. {\em
  Mathematical Biosciences} \textbf{179}, 21--55.

\bibitem{arefin2020social}
Arefin MR, Kabir KA, Jusup M, Ito H, Tanimoto J. 2020  Social efficiency
  deficit deciphers social dilemmas. {\em Scientific Reports} \textbf{10},
  16092.

\bibitem{wang2015universal}
Wang Z, Kokubo S, Jusup M, Tanimoto J. 2015  Universal scaling for the dilemma
  strength in evolutionary games. {\em Physics of Life Reviews} \textbf{14},
  1--30.

\bibitem{ito2018scaling}
Ito H, Tanimoto J. 2018  Scaling the phase-planes of social dilemma strengths
  shows game-class changes in the five rules governing the evolution of
  cooperation. {\em Royal Society Open Science} \textbf{5}, 181085.

\bibitem{tanimoto2007relationship}
Tanimoto J, Sagara H. 2007  Relationship between dilemma occurrence and the
  existence of a weakly dominant strategy in a two-player symmetric game. {\em
  BioSystems} \textbf{90}, 105--114.

\bibitem{arefin2021imitation}
Arefin MR, Tanimoto J. 2021  Imitation and aspiration dynamics bring different
  evolutionary outcomes in feedback-evolving games. {\em Proceedings of the
  Royal Society A} \textbf{477}, 20210240.

\bibitem{traulsen2023future}
Traulsen A, Glynatsi NE. 2023  The future of theoretical evolutionary game
  theory. {\em Philosophical Transactions of the Royal Society B: Biological
  Sciences} \textbf{378}, 20210508.

\bibitem{yc1}
Tavoni A. 2013  Building up cooperation. {\em Nature Climate Change}
  \textbf{3}, 782--783.

\bibitem{first2}
Wang X, Gu C, Quan J. 2019  Evolutionary game dynamics of the Wright-Fisher
  process with different selection intensities. {\em Journal of Theoretical
  Biology} \textbf{465}, 17--26.

\bibitem{second1}
Han TA, Tran-Thanh L. 2018  Cost-effective external interference for promoting
  the evolution of cooperation. {\em Scientific Reports} \textbf{8}, 15997.

\bibitem{second2}
Pereira LM, Anh HT. 2009  Evolution prospection. In {\em New Advances in
  Intelligent Decision Technologies: Results of the First KES International
  Symposium IDT 2009} pp. 51--63. Springer.

\bibitem{third1}
Kabir KA, Tanimoto J. 2020  Evolutionary game theory modelling to represent the
  behavioural dynamics of economic shutdowns and shield immunity in the
  COVID-19 pandemic. {\em Royal Society Open Science} \textbf{7}, 201095.

\bibitem{wang2020eco}
Wang X, Fu F. 2020  Eco-evolutionary dynamics with environmental feedback:
  Cooperation in a changing world. {\em Europhysics Letters} \textbf{132},
  10001.

\bibitem{ostrom1990governing}
Ostrom E. 1990 {\em Governing the commons: The evolution of institutions for
  collective action}.
Cambridge University Press.

\bibitem{leach1999environmental}
Leach M, Mearns R, Scoones I. 1999  Environmental entitlements: dynamics and
  institutions in community-based natural resource management. {\em World
  Development} \textbf{27}, 225--247.

\bibitem{hua2023evolution}
Hua S, Hui Z, Liu L. 2023  Evolution of conditional cooperation in
  collective-risk social dilemma with repeated group interactions. {\em
  Proceedings of the Royal Society B: Biological Sciences} \textbf{290},
  20230949.

\bibitem{sarkar2023managing}
Sarkar S. 2023  Managing ecological thresholds of a risky commons. {\em Royal
  Society Open Science} \textbf{10}, 230969.

\bibitem{hilbe2018evolution}
Hilbe C, {\v{S}}imsa {\v{S}}, Chatterjee K, Nowak MA. 2018  Evolution of
  cooperation in stochastic games. {\em Nature} \textbf{559}, 246--249.

\bibitem{wang2016effects}
Wang Y, Chen C, Araral E. 2016  The effects of migration on collective action
  in the commons: Evidence from rural China. {\em World Development}
  \textbf{88}, 79--93.

\bibitem{shu2023determinants}
Shu L, Fu F. 2023  Determinants of successful mitigation in coupled
  social-climate dynamics. {\em Proceedings of the Royal Society A}
  \textbf{479}, 20230679.

\bibitem{han2022institutional}
Han TA. 2022  Institutional incentives for the evolution of committed
  cooperation: ensuring participation is as important as enhancing compliance.
  {\em Journal of the Royal Society Interface} \textbf{19}, 20220036.

\bibitem{duong2021cost}
Duong MH, Han TA. 2021  Cost efficiency of institutional incentives for
  promoting cooperation in finite populations. {\em Proceedings of the Royal
  Society A} \textbf{477}, 20210568.

\bibitem{wang2019exploring}
Wang S, Chen X, Szolnoki A. 2019  Exploring optimal institutional incentives
  for public cooperation. {\em Communications in Nonlinear Science and
  Numerical Simulation} \textbf{79}, 104914.

\bibitem{cimpeanu2023does}
Cimpeanu T, Santos FC, Han TA. 2023  Does spending more always ensure higher
  cooperation? An analysis of institutional incentives on heterogeneous
  networks. {\em Dynamic Games and Applications} \textbf{13}, 1236--1255.

\bibitem{duong2023cost}
Duong MH, Durbac C, Han T. 2023  Cost optimisation of hybrid institutional
  incentives for promoting cooperation in finite populations. {\em Journal of
  Mathematical Biology} \textbf{87}, 77.

\end{thebibliography}
\end{document}